\documentclass[reqno,12pt]{amsart}
\newcommand{\CC}{{\mathbb C}}

\newcommand{\NN}{{\mathbb N}}

\def\bege{\begin{equation}} \def\ende{\end{equation}}
\def\begr{\begin{eqnarray}} \def\endr{\end{eqnarray}}
\usepackage{amsmath}
\usepackage{amssymb}
\usepackage{cite}

\usepackage{bbm}
\usepackage{amsmath}
\usepackage{txfonts}
\usepackage{stmaryrd}
\usepackage{amssymb}
\usepackage{amsfonts}
\usepackage{times}
\usepackage{mathrsfs}

\usepackage{amsthm}
\usepackage{enumerate}

\usepackage{framed}
\usepackage{lipsum}
\usepackage{color}
\usepackage{graphicx}

\newcommand{\TT}{{\mathbb T}}

\def\CC{ \mathbb{C}}
\newcommand{\DD}{{\mathbb D}}
\def\B{\mathcal{B}}

\def\D{\mathbb{D}}
\def\a{\alpha}

\def\t{\theta}

\def\begr{\begin{eqnarray}} \def\endr{\end{eqnarray}}

\def\msk{\medskip}
\def\ol{\overline}
\newtheorem{Lemma}{Lemma}[section]

\newtheorem{Proposition}[Lemma]{Proposition}
\newtheorem{Definition}[Lemma]{Definition}

\newcounter{other}            % Questions get letters
\begin{document}
	\title[]{On a conjecture about generalized integration operators on Hardy spaces}
	
	\author{Rong Yang and Songxiao Li$^\dagger$  }
	
	\address{Rong Yang
	\\ Department of Mathematics, Shantou University, 515063, Shantou, Guangdong,  P.R. China.\newline\indent
	Institute of Fundamental and Frontier Sciences, University of Electronic Science and Technology of China, 610054, Chengdu, Sichuan, P.R. China.}
	\email{yangrong071428@163.com  }
	
	\address{Songxiao Li\\ Department of Mathematics, Shantou University, 515063, Shantou,  Guangdong, P.R. China. } \email{jyulsx@163.com}

%\begin{figure}[hp]
% \begin{center}
%\centering\includegraphics[scale=0.75]{email.png}
% \end{center} 
%\end{figure}

	
	\begin{abstract} 
	A conjecture posed by Chalmoukis in 2020 states that if $T_{g,a}:H^p\to H^q(0<q<p<\infty)$ is bounded, then $g$ must be in $H^{\frac{pq}{p-q}}$.
	In this article, we provide a positive answer to the  aforementioned conjecture. We also consider the compactness of $T_{g,a}:H^p\to H^q(0<q<p<\infty)$.

	\thanks{$\dagger$ Corresponding author.}
	\end{abstract}
	
	\maketitle
	
\section{Introduction}
Let $\mathbb{D}$ denote the unit disc on the complex plane $\mathbb{C}$, $\mathbb{T}$ its boundary, and $H(\mathbb{D})$ be the set of
all analytic functions in $\mathbb{D}$. We denote by $\NN$ the set of positive integers.
The Hardy space $H^p(0<p<\infty)$ consists of all $f \in H(\mathbb{D})$ for which
$$
\|f\|_{H^p}^p=\sup _{0\leq r<1} \int_0^{2 \pi}|f(r e^{i \theta})|^p \frac{d\t}{2 \pi}<\infty.
$$
	
Let $g \in H(\mathbb{D})$, $k\in\NN\cup\{0\}$  and $n\in\NN$ such that $0 \leq k<n$. The generalized integration operator $T_g^{n, k}$ is defined by
$$
T_g^{n, k} f(z)=I^n\left(f^{(k)}(z) g^{(n-k)}(z)\right), \quad f \in H(\mathbb{D}).
$$
Here,  $I^n$ is the $n$-th iteration of the integration operator $I f(z)=\int_0^z f(t)dt$.
The operator $T_g^{n, k}$ was first introduced by Chalmoukis \cite{ch}. In particular, when $n=1$ and $k=0$, we have
$$
T_g^{1,0} f(z)=\int_0^z f(\eta) g^{\prime}(\eta) d \eta=T_g f(z) .
$$	

In 1977, the Volterra integration operator $T_g$ was first introduced and stuided by Pommerenke in \cite{p}. It was proved that $T_g$ is bounded on $H^2$ if and only if $g\in BMOA$. In 1995, Aleman and Siskakis proved that $T_g$ is bounded  on $H^p$($p\ge1$) if and only if $g\in BMOA$ in \cite{as1}.
In \cite{as2}, they proved that $T_g$ is bounded  on $A^p$ if and only if $g\in\B$, the Bloch space. 
See \cite{ag1, mppw,  llw} and the references therein for more study of the operator $T_g$.
In \cite{ch},  Chalmoukis studied the boundedness of  the  operator $T_{g}^{n,k}$ between different Hardy spaces. For example, he showed that $T_g^{n,0}: H^p\to H^p$ is bounded if and only if $g\in BMOA$, $T_g^{n,k}: H^p\to H^p$ is bounded  if and only if $g\in\B$ when $k\geq 1$, while $T_g^{n,k}:H^p\to H^q$ is bounded  if and only if
$\sup\limits_{z\in\D}(1-|z|^2)^{\frac{1}{q}-\frac{1}{p}+n-k}|g^{(n-k)}(z)|<\infty$ when $0<p<q<\infty$. 
From the results just mentioned, we see that the operator $T_g^{n,k}$ has a very different behavior from that of $T_g$.

Let $g\in H(\D)$ and $a=(a_1,a_2,\cdots,a_{n-1})\in\CC^{n-1}$. The generalized integration operator $T_{g,a}$ is defined by (see \cite{ch})
\begin{align*}
T_{g,a}=T_g^{n,0}+\sum_{k=1}^{n-1}a_kT_g^{n,k}.
\end{align*}
Chalmoukis proposed a conjecture in \cite{ch}.
Is it true that if $T_{g,a}:H^p\to H^q$, $0<q<p<\infty$ is bounded, then $g$ must be in $H^s$, $\frac{1}{s}=\frac{1}{q}-\frac{1}{p}?$

Due to the relationship between $T_{g,a}$ and $T_g^{n,k}$, to characterize the boundedness of $T_{g,a}$, it suffices to characterize the boundedness of $T_g^{n,k}$ for any $n\in\NN$  and $k\in\NN\cup\{0\}$ such that $0\le  k< n$. In this paper, we give a positive answer to this conjecture. More specifically, the main result of this paper can be stated as follows. \msk

\noindent{\bf Main Theorem.} Let $0<q<p<\infty$, $g\in H(\D)$, $n\in\NN$  and $k\in\NN\cup\{0\}$ such that $0\le  k< n$. Then the following statements are equivalent:
	
$(i)$ The operator $T_g^{n,k}:H^p\to H^q$ is bounded.
	
$(ii)$ The operator $T_g^{n,k}:H^p\to H^q$ is compact.
	
$(iii)$  $g\in H^s$, where $\frac{1}{s}=\frac{1}{q}-\frac{1}{p}$. \msk
 
This paper is organized as follows. In Section 2, we mainly state some definitions and preliminary results. In Section 3, we provide the proof of the main theorem.
	
	In this paper, the symbol $A\lesssim B$ means that there exists a positive finite constant $C$ such that $A\le CB$. The value of $C$ may change from one occurrence to another. The expression $A\asymp B$ means that both $A\lesssim B$ and $B\lesssim A$ hold.\msk

\section{Preliminaries}

In this section, we state some notations, lemmas and propositions, which will be used in the proof of main results in this paper.

Let $\beta(z, w)$ be the hyperbolic metric on $\mathbb{D}$. 
%that is, 
%	$$
%	\b(z,w)=\frac{1}{2}\log\frac{1+\rho(z,w)}{1-\rho(z,w)},\quad z,w\in\D.
%	$$
If $\beta(z, w)<r$ for any $z, w \in \mathbb{D}$, then $ 1-|w| \asymp 1-|z|$.
Let $D(z, r)=\{w \in \mathbb{D}: \beta(z, w)<r\}$ be the hyperbolic disc of radius $r>0$ centered at $z \in \mathbb{D}$. 
Let $\eta \in \mathbb{T}$ and $\zeta>\frac{1}{2}$. The non-tangential region $\Gamma_\zeta(\eta)$ is defined by
$$
\Gamma(\eta)=\Gamma_\zeta(\eta) =\{z \in \mathbb{D}:|z-\eta|<\zeta(1-|z|^2) \}.
$$
The following lemma is very important in this paper and will be used frequently.

\begin{Lemma}\cite[Lemma 2.3]{w1}\label{2.9}
	Let $\zeta>1, r \geq 0$ and $\eta \in \mathbb{T}$. If $\zeta_+=(\zeta+1) e^{2 r}-1$, then
	\begin{align}\label{+}
		D(z, r) \subset \Gamma_{\zeta_+}(\eta)
	\end{align}
	for all $z \in \Gamma_{\zeta}(\eta)$.

\end{Lemma}
In the rest of this paper, we will denote  $\Gamma_{\zeta_{+}}(\eta)$ by $\widetilde{\Gamma}(\eta)$ in (\ref{+}) when we need not mention the aperture $\zeta_+$.

\begin{Definition}
	Let $0<p, q<\infty$ and $\a>-2$. The tent space $T_p^q(\alpha)$ consists of all measurable functions $f$ on $\mathbb{D}$ with
	\begin{align*}
		\|f\|_{T_p^q(\alpha)}=\left(\int_\mathbb{T}\left(\int_{\Gamma(\eta)}|f(z)|^p(1-|z|^2)^{\a} d A(z)\right)^{\frac{q}{p}}|d \eta|\right)^{\frac{1}{q}}<\infty .
	\end{align*}
	In particular, for $\a=0$, we write $T_p^q$ instead of $T_p^q(\a)$.
	
	For $p=\infty$ and $0<q<\infty$, the tent space $T_{\infty}^q(\alpha)$ consists of all measurable functions $f$ on $\mathbb{D}$ such that
	\begin{align*}%\label{2.222}
	\|f\|_{T_\infty^q(\alpha)}=\left(\int_\mathbb{T} \left(\operatorname{esssup}_{z \in \Gamma(\eta)}|f(z)|\right)^q|d \eta|\right)^{\frac{1}{q}}<\infty.
	\end{align*}

	For $q=\infty$ and $0<p<\infty$,
	the tent space $T_p^{\infty}(\alpha)$ consists of all measurable functions $f$ on $\mathbb{D}$ for which
	$$
	\|f\|_{T_p^{\infty}(\alpha)}=\operatorname{esssup} _{\eta \in \mathbb{T}}\left(\sup _{u \in \Gamma(\eta)} \frac{1}{1-|u|^2} \int_{S(u)}|f(z)|^p(1-|z|^2)^{\a+1} d A(z)\right)^{\frac{1}{p}}<\infty,
	$$
	where
	$$
	S(r e^{i \theta})=\left\{\lambda e^{i t}: |t-\theta| \leq \frac{1-r}{2},1-\lambda \leq 1-r\right\}
	$$
	for $r e^{i \theta} \in \mathbb{D} \backslash\{0\}$ and $S(0)=\mathbb{D}$.
\end{Definition}

Tent spaces were initially introduced by Coifman, Meyer and Stein  in \cite{cms} to study problems in harmonic analysis. They provided a  general framework for questions regarding some important spaces such as Hardy spaces and Bergman spaces. The aperture $\zeta  $ of the non-tangential region $\Gamma_\zeta(\eta)$   is suppressed in the above definition since   any two apertures generate the same function space with equivalent quasinorms.	Denote $T_p^q(\alpha) \cap H(\mathbb{D})$ by $A T_p^q(\alpha)$.
Notice that the definition of $T^q_\infty(\a)=T^q_\infty$ is independent of $\a$. 
By the non-tangential maximal function characterization of the Hardy space, we see that  $AT^q_\infty=H^q$(see \cite{z2}). 
 	From \cite[Theorem G]{pa}, we have the following lemma.
	\begin{Lemma}\label{H}
	Let $0<s<\infty$ and $n\in\NN$.  If $f\in H(\D)$ and $f(0)=0$,  then
	$$
	\|f\|_{H^s}^s\asymp\int_{\TT}\left(\int_{\Gamma(\eta)}|f^{(n)}(z)|^2(1-|z|^2)^{2n-2}dA(z)   \right)^{\frac{s}{2}}|d\eta|.
	$$
	\end{Lemma}
 
The sequence $Z=\{a_j\}$ is a separated sequence if there is a constant $\tau>0$ with $\beta(a_j, a_j) \geq \tau$ for $k \neq j$. 
For $r>\kappa>0$, the sequence $Z=\{a_j\}$ is said to be an $(r, \kappa)$-lattice if 
$\mathbb{D}=\bigcup_k D(a_j, r)$ and
the sets $D(a_j, \kappa)$ are pairwise disjoint.
It is obvious that any $(r, \kappa)$-lattice is a separated sequence.

\begin{Definition}
Let $0<p, q<\infty$ and $Z=\{a_j\}$ be an $(r, \kappa)$-lattice. The tent space $T_p^q(Z)$ consists of all $\{x_j\}$ such that
\begin{align*}%\label{(2.4)}
\|\{x_j\}\|_{T_p^q(Z)}=\left(\int_\mathbb{T}\left(\sum_{a_j \in \Gamma(\eta)}|x_j|^p\right)^{\frac{q}{p}}|d \eta|\right)^{\frac{1}{q}}<\infty.
\end{align*}
For $p=\infty$ and $0<q<\infty$, the tent space $T_\infty^q(Z)$ consists of all $\{x_j\}$ with
$$
\|\{x_j\}\|_{T_\infty^q(Z)}=\left(\int_\mathbb{T}\sup _{a_j \in \Gamma(\eta)}|x_j|^q|d \eta|\right)^{\frac{1}{q}}<\infty.
$$
For $q=\infty$ and $0<p<\infty$, the tent space $T^\infty_p(Z)$ consists of all $\{x_j\}$ for which
$$
\|\{x_j\}\|_{T_p^{\infty}(Z)}=\operatorname{esssup} _{\eta \in \mathbb{T}}\left(\sup _{u \in \Gamma(\eta)} \frac{1}{1-|u|^2} \sum_{a_j \in S(u)}|x_j|^p(1-|a_j|^2)\right)^{\frac{1}{p}}<\infty.
$$
\end{Definition}
	
Analogously, the aperture $\zeta  $ of the non-tangential region $\Gamma_\zeta(\eta)$   is suppressed in the above definition since any two apertures generate the tent spaces of sequences  with equivalent quasinorms.

\begin{Lemma}\cite[Lemma 3]{mppw}\label{za}
Let $0<p<\infty$ and $\a\ge0$. There exists $r_0\in(0,1)$ so that if $0<r<r_0$ and
$Z=\{a_j\}$ is an $(r,\kappa)$-lattice, then
$$
\int_\mathbb{T} \sup_{z \in \Gamma(\eta)}|f(z)|^p (1-|z|^2)^\a|d \eta|\lesssim \int_\mathbb{T} \sup_{a_j \in \Gamma(\eta)}|f(a_j)|^p (1-|a_j|^2)^\a|d \eta|, 
$$
whenever $f$ is analytic on $\D$ such that the left-hand-side is finite.
\end{Lemma}

The following lemma is very important, and we can find it in \cite{pj}. % [Proposition 7.2]
	
\begin{Proposition}\label{7.2}
Let $0<p<\infty$, $b>\max \{1, \frac{2}{p}\}$, and
$Z=\{a_j\}$ be an $(r, \kappa)$-lattice. Then the operator $T_Z: T_2^p(Z) \rightarrow  H^p$ is bounded, where,
$$
T_Z(\{x_j\})(z)=\sum_{j=1}^{\infty} x_j \frac{(1-|a_j|)^{b}}{(1-\overline{a_j} z)^{b}}.
$$
\end{Proposition}

We will use the following proposition, which is a factorization  of tent spaces  of sequences.

\begin{Proposition}\cite[Proposition 2.16]{ag1}\label{2.16}
Let $0<p, q<\infty$ and $Z=\{a_j\}$ be an $(r, \kappa)$-lattice. If $p \leq p_1, p_2 \leq\infty, q \leq q_1$, $q_2 \leq\infty$ and satisfy $\frac{1}{p}=\frac{1}{p_1}+\frac{1}{p_2}$ and $\frac{1}{q}=\frac{1}{q_1}+\frac{1}{q_1}$. Then
$$
T_p^q(Z)=T_{p_1}^{q_1}(Z) \cdot T_{p_2}^{q_2}(Z).
$$
\end{Proposition}

The following proposition play a key role in the proof of our main theorem. 
From now on, $p'$ is the conjugate exponent of $p$ with $\frac{1}{p}+\frac{1}{p'}=1$. In this paper, we agree that $1'=\infty$ and $\infty'=1$.

\begin{Proposition}\cite[Lemma 6]{ar}\label{2.18}
Let $1 \leq p<\infty, 1<q<\infty$ and $Z=\{a_j\}$ be an $(r, \kappa)$-lattice . Then $\left(T_p^q(Z)\right)^* \cong T_{p^{\prime}}^{q^{\prime}}(Z)$. The isomorphism between $\left(T_p^q(Z)\right)^*$ and $T_{p^{\prime}}^{q^{\prime}}(Z)$ is given by the operator
$
\{d_j\} \mapsto\langle\cdot,\{d_j\}\rangle,
$
where $\langle\cdot,\{d_j\}\rangle$ is defined by
$$
\langle\{c_j\},\{d_j\}\rangle=\sum_j c_j d_j(1-|a_j|^2), \quad\{c_j\} \in T_p^q(Z).
$$
In fact, $$\|\{d_j\}\|_{T_{p^{\prime}}^{q^{\prime}}(Z)} \asymp \sup \left\{\left|\sum_j c_j d_j(1-|a_j|^2)\right|:\|\{c_j\}\|_{T_p^q(Z)}=1\right\}.$$ 
\end{Proposition}

\begin{Proposition}\cite[Proposition 2.7]{ag1}\label{2.7}
Let $Z=\{a_j\}$ be an $(r, \kappa)$-lattice and $L \geq 1, R>0$. There is a positive integer $N=$ $N(L, R, Z)$ such that for each point $z \in \mathbb{D}$ there are at most $N$ hyperbolic discs $D(a_j, L r)$ satisfying $D(z, R) \cap D(a_j, L r) \neq \emptyset$.
\end{Proposition}\msk

\section{Proof of the main theorem}
\noindent {\bf Proof.} $(ii)\Rightarrow (i).$ It is obvious.

$(i)\Rightarrow (iii).$
Let  $Z=\{a_j\}$ be an $(r, \kappa)$-lattice and $\{x_j\}\in T_2^p(Z)$. Let $r_j:[0,1] \rightarrow\{-1,1\}$ be the Radermacher functions. For $L>\max \{1, \frac{2}{p}\}$, set
$$
F_u(z)=\sum_{j=1}^\infty x_j r_j(u) \frac{(1-|a_j|^2)^{L}}{(1-\overline{a_j} z)^{L}}, \quad z \in \mathbb{D} .
$$
Using Proposition $\ref{7.2}$, we get $
\|F_u\|_{H^p}\lesssim\|\{x_j\}\|_{T_2^p(Z)}.$
By the assumption and Lemma \ref{H}, we get
\begin{align*} &\int_\mathbb{T}\left(\int_{\Gamma(\eta)}|(T_g^{n,k}F_u)^{(n)}(z)|^2(1-|z|^2)^{2n-2} dA(z)\right)^{\frac{q}{2}}|d \eta|\\
\lesssim& \|T_g^{n,k}\|^q \|F_u\|_{H^p}^q
\lesssim \|T_g^{n,k}\|^q \|\{x_j\}\|_{T_2^p(Z)}^q.
\end{align*}
Integrating both sides with respect to $u$ from $0$ to $1$,
using Fubini's theorem and Kahane's inequality (see \cite[Theorem 2.1]{ka}), we have
\begin{align*}
&\int_\mathbb{T}\left(\int_{0}^{1}\int_{\Gamma(\eta)}
|F_u^{(k)}(z)|^2|g^{(n-k)}(z)|^2
(1-|z|^2)^{2n-2} dA(z) du \right)^{\frac{q}{2}} |d \eta|\\
\lesssim&\int_\mathbb{T}\int_{0}^{1}\left(\int_{\Gamma(\eta)}|F_u^{(k)}(z)|^2|g^{(n-k)}(z)|^2(1-|z|^2)^{2n-2} dA(z)\right)^{\frac{q}{2}} du |d \eta|\\ 
\lesssim& \|T_g^{n,k}\|^q \|\{x_j\}\|_{T_2^p(Z)}^q.
\end{align*}
Applying Khinchine's inequality (see \cite{l3}), we obtain
%\begin{align*}
%&\int_\mathbb{T}\left(\int_{\Gamma(\eta)}\left(\sum_j|x_j|^2 |a_j|^{2k}\frac{(1-|a_j|^2)^{2 L}}{|1-\overline{a_j} z|^{2 L+2k}}\right) |g^{n-k}(z)|^2
%(1-|z|^2)^{2n-2} dA(z)\right)^{\frac{q}{2}}|d \eta| \\
%\lesssim& \|T_g^{n,k}\|^q\|\{x_j\}\|_{T_2^p(Z)}^q.
%\end{align*}
%Hence, 
\begin{align*}
&\int_\mathbb{T}\left(\int_{\Gamma(\eta)}\left(\sum_j|x_j|^2 \frac{(1-|a_j|^2)^{2 L}}{|1-\overline{a_j} z|^{2 L+2k}}\right) |g^{n-k}(z)|^2
(1-|z|^2)^{2n-2} dA(z)\right)^{\frac{q}{2}}|d \eta| \\
\lesssim& \|T_g^{n,k}\|^q\|\{x_j\}\|_{T_2^p(Z)}^q.
\end{align*}
Let $w_j=|g^{(n-k)}(a_j)|(1-|a_j|^2)^{n-k}$. 
From Lemma \ref{2.9} we get  $\cup_{z\in\Gamma(\eta)}D(z,4r)\subset\widetilde{\Gamma}(\eta)$. Using Proposition \ref{2.7}, Fubini's theorem and H\"older's inequality, we get
\begin{align*}
&\int_\mathbb{T}\left(\sum_{a_j \in \Gamma(\eta)}|x_j|^2 |w_j|^2\right)^{\frac{q}{2}}|d \eta| \\ 
=&\int_\mathbb{T}\left(\sum_{a_j \in \Gamma(\eta)}|x_j|^2 |g^{(n-k)}(a_j)|^2(1-|a_j|^2)^{2n-2k}\right)^{\frac{q}{2}}|d \eta| \\ 
\lesssim&\int_\mathbb{T}\left(\sum_{a_j \in \Gamma(\eta)}|x_j|^2 \int_{D(a_j, 4r)}|g^{(n-k)}(z)|^2 \frac{(1-|z|^2)^{2n-2}(1-|a_j|^2)^{2L}}
{|1-\ol{a_j}z|^{2L+2k}}dA(z)\right)^{\frac{q}{2}}|d \eta| \\
\lesssim&\int_{\mathbb{T}}\left(\int_{\widetilde{\Gamma}(\eta)}\sum_j|x_j|^2 
\frac{(1-|a_j|^2)^{2L}}{|1-\ol{a_j}z|^{2L+2k}}\chi_{D(a_j,4r)}(z)|g^{(n-k)}(z)|^2(1-|z|^2)^{2n-2}dA(z)\right)^{\frac{q}{2}}|d\eta|\\
\lesssim&\int_{\mathbb{T}}\left(\int_{\widetilde{\Gamma}(\eta)}\left(\sum_j|x_j|^2 
\frac{(1-|a_j|^2)^{2L}}{|1-\ol{a_j}z|^{2L+2k}}\right)|g^{(n-k)}(z)|^2(1-|z|^2)^{2n-2}dA(z)\right)^{\frac{q}{2}}|d \eta|\\
\lesssim&\|T_g^{n,k}\|^q \|\{x_j\}\|_{T_2^p(Z) }^q.
\end{align*}
Let $b>\max \left\{\frac{1}{q}, \frac{1}{2}\right\}$. For any $\{l_j\} \in T_{\frac{2b}{2b-1}}^{\frac{bq}{bq-1}}(Z)$,
using Fubini's theorem and H\"older's inequality twice, we obtain
\begin{equation}\label{4.5}
\begin{aligned}
&\sum_j |l_j x_j^{\frac{1}{b}} w_j^{\frac{1}{b}}|(1-|a_j|^2) \\
\asymp& \int_{\mathbb{T}}\left(\sum_{a_j \in \Gamma(\eta)} |l_j| |x_j|^{\frac{1}{b}} |w_j|^{\frac{1}{b}}\right)|d \eta| \\
\lesssim& \int_{\mathbb{T}}\left(\sum_{a_j \in \Gamma(\eta)}|x_j|^2  |w_j|^2 \right)^{\frac{1}{2b}}\left(\sum_{a_j \in \Gamma(\eta)} |l_j|^{\frac{2b}{2b-1}}\right)^{\frac{2b-1}{2b}}|d \eta| \\
\leq&\left(\int_{\mathbb{T}}\left(\sum_{a_j \in \Gamma(\eta)}|x_j|^2 |w_j|^2\right)^{\frac{q}{2}}|d \eta|\right)^{\frac{1}{bq}}
\left( \int_{\mathbb{T}} \left(  \sum_{a_j \in \Gamma(\eta)} |l_j|^{\frac{2b}{2b-1}} \right)^\frac{q(2b-1)}{2(bq-1)}   |d\eta| \right)^{\frac{bq-1}{bq}}\\
\lesssim&\|T_g^{n,k}\|^{\frac{1}{b}} \|\{x_j\}\|_{T_2^p(Z)}^{\frac{1}{b}}\|\{l_j\}\|_{T_{\frac{2b}{2b-1}}^{\frac{bq}{bq-1}}(Z)}.
\end{aligned}
\end{equation}
Applying Proposition \ref{2.16} we get
$$
\{l_j x_j^\frac{1}{b}\} \in T_{1}^{\frac{q b s}{pbq-p+q}}(Z)=T_{\frac{2b}{2b-1}}^{\frac{bq}{b q-1}}(Z) \cdot T_{2b}^{pb}(Z).
$$
By Proposition \ref{2.18},
we obtain 
$$\left(T_{1}^{\frac{pbq}{pbq-p+q}}(Z)\right)^*=T_{\infty}^{\frac{pbq}{p-q}}(Z).  
$$
By $(\ref{4.5})$ we have
\begin{align*}
&\left\|\left\{|w_j|^{\frac{1}{b}}\right\}\right\|_{T_{\infty}^{\frac{pbq}{p-q}}(Z)}=\left\|\left\{|w_j|\right\}\right\|_{T_{\infty}^{\frac{pq}{p-q}}(Z)}^\frac{1}{b}\\
=&\left(\int_{\mathbb{T}}\left(\sup _{a_j \in \Gamma(\eta)} |g^{(n-k)}(a_j)|(1-|a_j|^2)^{n-k}\right)^{\frac{pq}{p-q}}|d \eta|\right)^{\frac{p-q}{pbq}}<\infty .
\end{align*}
Thus, we have $|g^{(n-k)}(a_j)|(1-|a_j|^2)^{n-k}\in T_\infty^{\frac{pq}{p-q}}(Z)=T^s_\infty(Z)$.
By Lemma \ref{za}, we obtain
$g^{(n-k)}(\cdot)(1-|\cdot|^2)^{n-k}\in AT_\infty^s$, i.e., $g\in H^s$.

$(iii)\Rightarrow (ii).$
Assume that
$g\in H^s$. Then $g^{(n-k)}(\cdot)(1-|\cdot|^2)^{n-k}\in AT_\infty^{s}=H^s$.  
Using Lemma \ref{H} and H\"older's inequality, we have
\begin{equation}\label{bb}
\begin{aligned}
&\|T_g^{n,k}f\|_{H^q}^q\\
\asymp&\int_{\mathbb{T}}\left( \int_{\Gamma(\eta)}|(T_g^{n,k}f)^{(n)}(z)|^2(1-|z|^2)^{2n-2}dA(z)  \right)^{\frac{q}{2}}|d\eta|\\
=&\int_{\mathbb{T}}\left( \int_{\Gamma(\eta)}|f^{(k)}(z)|^2|g^{(n-k)}(z)|^2(1-|z|^2)^{2n-2}dA(z)  \right)^{\frac{q}{2}}|d\eta|\\
\lesssim& \int_{\mathbb{T}} \left(\sup_{z\in \Gamma(\eta)}|g^{(n-k)}(z)|^{2}(1-|z|^2)^{2n-2k}  \right)^{\frac{q}{2}}\cdot \left(   \int_{\Gamma(\eta)}|f^{(k)}(z)|^2(1-|z|^2)^{2k-2}dA(z) \right)^{\frac{q}{2}}
|d\eta| \\
\lesssim&
\left(\int_{\mathbb{T}} \left(\sup_{z\in \Gamma(\eta)}|g^{(n-k)}(z)|(1-|z|^2)^{n-k}  \right)^{\frac{pq}{p-q}}|d\eta|  \right)^{\frac{p-q}{p}}\\
&\cdot
\left(  \int_{\mathbb{T}} \left(   \int_{\Gamma(\eta)}|f^{(k)}(z)|^2(1-|z|^2)^{2k-2}dA(z) \right)^{\frac{p}{2}} |d\eta|\right)^{\frac{q}{p}}\\
\lesssim& \|g^{(n-k)}\|^q_{AT_\infty^{s}} \left(  \int_{\mathbb{T}} \left(   \int_{\Gamma(\eta)}|f^{(k)}(z)|^2(1-|z|^2)^{2k-2}dA(z) \right)^{\frac{p}{2}} |d\eta|\right)^{\frac{q}{p}}\\
\lesssim& \|g^{(n-k)}\|^q_{H^s}\|f\|^q_{H^p}.
\end{aligned}
\end{equation}
Thus, $T_g^{n,k}:H^p\to H^q$ is bounded.

Consider the sequence $\{f_j\}_{j=1}^{\infty} \subset H^p$ such that $\sup_j \|f_j\|_{H^p} < \infty$. This ensures that $\{f_j\}$ is uniformly bounded on compact subsets of  $\mathbb{D}$. By Montel's theorem, it forms a normal family. Hence, we can extract a subsequence $\{f_{n_j}\}_{j=1}^{\infty}$ which converges uniformly on compact subsets of $\mathbb{D}$ to $f\in H(\D)$. According to Fatou's Lemma, $f \in H^p$. Define $h_j = f_{n_j} - f$. Then $h_j\in H^p$. Now, we only need to prove that 
$$\lim_{j \to \infty} \|T_g^{n,k} h_j\|_{H^q} = 0,$$
which implies that $T_g^{n,k}: H^p \to H^q$ is compact.
From the fact that
$g^{(n-k)}(\cdot)(1-|\cdot|^2)^{n-k}\in AT_\infty^{s}=H^s$,
%we have
%$$
%\lim_{R\to 1^{-}}\int_{\TT}\left( \sup_{{z}\in\Gamma(\eta)\backslash \ol{D(0,R)}} |g^{(n-k)}(z)|(1-|z|^2)^{n-k} \right)^{\frac{pq}{p-q}}|d\eta|=0,
%$$
%which implies that 
for any $\epsilon>0$, there exists $R_0\in(0,1)$ such that
$$
\left(\int_{\TT}\left( \sup_{z\in\Gamma(\eta)\backslash \ol{D(0,R)}} |g^{(n-k)}(z)|(1-|z|^2)^{n-k} \right)^{\frac{pq}{p-q}}|d\eta|\right)^{\frac{p-q}{pq}}<\epsilon
$$
when $R\ge R_0$.
Choose $j_0$ such that $\sup_{j\ge j_0,|z|\le R_0}|h_j^{(k)}(z)|<\epsilon$. Using $(\ref{bb})$,  we get
\begin{align*}
 \|T_g^{n,k}h_j\|_{H^q}^q\lesssim&\int_{\mathbb{T}}\left( \int_{\Gamma(\eta)\cap\{|z|\le R_0\}}|h_j^{(k)}(z)|^2|g^{(n-k)}(z)|^2(1-|z|^2)^{2n-2}dA(z)  \right)^{\frac{q}{2}}|d\eta|\\
&+\int_{\mathbb{T}}\left( \int_{\Gamma(\eta)\backslash\ol{D(0,R_0)}}|h_j^{(k)}(z)|^2|g^{(n-k)}(z)|^2(1-|z|^2)^{2n-2}dA(z)  \right)^{\frac{q}{2}}|d\eta|\\
\lesssim&\epsilon^q+ \|h_j\|_{H^p}^q  \left(\int_{\TT}  \left(\sup_{z\in \Gamma(\eta)\backslash\ol{D(0,R_0)}}|g^{(n-k)}(z)|(1-|z|^2)^{n-k}  \right)^{\frac{pq}{p-q}}|d\eta| \right)^{\frac{p-q}{p}} \\
\lesssim&\epsilon^q
\end{align*}
when $j\ge j_0$.
Thus, $T_g^{n,k}:H^p\to H^q$ is compact. The proof is complete. \msk

	\noindent{\bf Data Availability}  No data was used to support this study.\msk
	
	\noindent{\bf Conflicts of Interest}  The authors  declare that they have no conflicts of interest.\msk

\noindent	{\bf Acknowledgements}  	The corresponding author was supported by  NNSF of China (No. 12371131), STU Scientific Research Initiation Grant(No. NTF23004), 
  Li Ka-Shing Foundation (Grant no. 2024LKSFG06), Guangdong Basic and Applied Basic Research Foundation (No. 2023A1515010614).

 % \begin{center}
%\centering\includegraphics[scale=0.4]{email.png}
% \end{center} 
 \begin{figure}[hp]
 \begin{center}
\centering\includegraphics[scale=0.75]{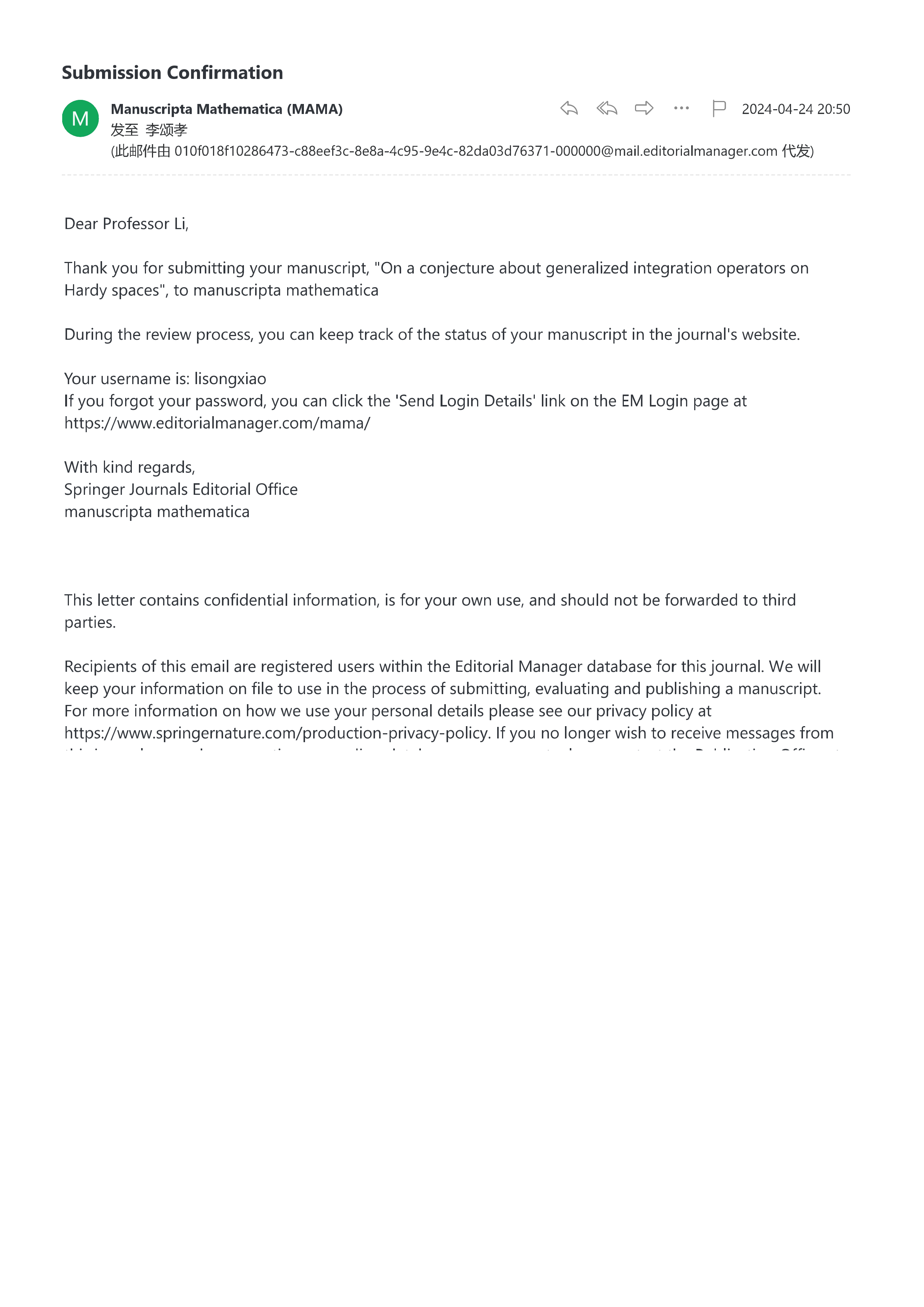}
 \end{center} 
\end{figure}

\end{document}